\documentclass{amsart}
\usepackage{amsfonts,amsmath,amssymb}
\newcommand{\vol}{\text{vol}}
\newcommand{\card}{\text{card }}
\newcommand{\proved}{\begin{flushright} \qed \end{flushright}}
\newtheorem{thm}{Theorem}
\newtheorem{prop}{Proposition}
\newtheorem{ex}{Example}
\begin{document}



\title{Kergin approximation in Banach spaces}


\author{Scott Simon}

\address{Mathematics Department, Stony Brook University, Stony Brook NY, 11794-3651, USA}

\newpage

\begin{abstract}
We explore the convergence of Kergin interpolation polynomials of holomorphic functions in Banach spaces, which need not be of bounded type.  We also investigate a case where the Kergin series diverges.
\end{abstract}
\maketitle


Kergin interpolation is a generalization of both Lagrange interpolation in the one dimensional case, and the Taylor polynomial in the case where all interpolation points coincide.  In several variables, interpolation polynomials are not unique.  However, Kergin \cite{Kergin} proved that interpolation polynomials enjoying natural properties exist and are unique:
\begin{thm}[Kergin]
Let $N \in \mathbb{N}^+, K \in \mathbb{N}$, and $x_0, \ldots, x_K \in \mathbb{R}^N$, not necessarily distinct.  There is a unique $\chi : \mathcal{C}^K (\mathbb{R}^N) \rightarrow P^K (\mathbb{R}^N)$ satisfying:
\begin{enumerate}
\item $\chi$ is linear.
\item For every $f \in \mathcal C^K (\mathbb{R}^N)$, every $q \in Q^k$ in $\mathbb{R}^N$, where $k \in \{0, \ldots, K\}$, and every $J \subset \{0, \ldots, K\}$ with $\card J= k+1$, there exists $x \in [x_j]_{j \in J}$ such that $q(\chi(f) - f)(x)=0$.
\end{enumerate}
\end{thm}
Here, $\mathcal{C}^K (\mathbb{R}^N)$ is the set of functions with $K$ continuous derivatives, $Q^k$ is the set of constant coefficient linear partial differential operators of order $k$, $P^K (\mathbb{R}^N)$ is the set of polynomials of degree $K$, and $[x_j]_{j \in J}$ is the convex hull of $\{x_j\}_{j \in J}.$  It fell to Micchelli \cite{Micchelli} and Milman \cite{MicchelliMilman} to discover a formula for these polynomials.  This formula also extends to the Banach space case, see \cite{Filipsson, Petersson}.  In this case, the potential unboundedness of continuous functions, even on bounded sets bounded away from the boundary of the domain, presents new difficulties in proving convergence results.  Filipsson \cite{Filipsson} proved a convergence result for holomorphic functions bounded on a ball.

We give the formula for the Kergin polynomial below.
Let $X, Y$ be complex Banach spaces, $U \subset X$ open and $f:U \rightarrow Y$.  Define
$d^0 f=f$ and $d^{k+1}f:U \times X^{k+1} \rightarrow Y$,
\[
d^{k+1}f(x;\xi_1, \ldots, \xi_{k+1})=\lim_{t \rightarrow 0} \frac{1}{t} d^k f(x+t\xi_{k+1};\xi_1, \ldots, \xi_k),
\]
if this limit exists.
This is just the $k+1^{st}$ iteration of the directional derivative of $f$, see, e.g., \cite{Lempert}.
Let $p_0, \ldots p_n, \in X.$  Suppose $f$ is an $n$-times differentiable function on the convex hull of $p_0, \ldots, p_n$. The Kergin polynomial of $f$ of degree $n$ is the sum
\begin{equation}\label{eq:kergindef}
f(p_0)+\sum_{k=1}^{n} \int_{S_k} d^k f(s_0 p_0+ \ldots +s_k p_k;x-p_0, \ldots, x-p_{k-1})ds_1 \ldots ds_k,
\end{equation}
where
\[
S_k=\{(s_1, \ldots, s_k) \in \mathbb{R}^k:s_j \ge 0, \sum_{j=1}^k s_j \le 1 \}
\]
is the standard $k$-simplex, and
\begin{equation}\label{eq:s_0}
s_0=1-s_1-\ldots - s_k.
\end{equation}
This is a Bochner, or vector valued, integral, see, e.g., \cite{Mujica}.  In the case where $Y=\mathbb{C},$ this is just the
usual Lebesgue (or Riemann) integration.

L. Filipsson observed that Micchelli's error formula for the degree $(k-1)$  Kergin polynomial carries over to Banach spaces:
\begin{equation}\label{eq:error}
\int_{S_k} d^k f(s_0 p_0+ \ldots +s_{k-1} p_{k-1} +s_k x;x-p_0,\ldots ,x-p_{k-1}) ds_1 \ldots ds_k,
\end{equation}
where $s_0$ is as in \eqref{eq:s_0}.
Given an infinite sequence $p_0, p_1, p_2, \ldots,$ we define the infinite Kergin series by replacing $n$ in \eqref{eq:kergindef} by $\infty$.  Under some circumstances, this series will approximate the given function.  That is the primary subject of this paper.  First, we will need the following proposition.
\begin{prop}\label{pr:cpct}
Let $X$ be a complex Banach space and let $K \subset X$ be compact.  Then
\begin{enumerate}
\item the convex hull of $K$ is compact, and
\item the balanced hull of $K$ is compact.
\end{enumerate}
\end{prop}
This can be found, for example, in \cite{Bourbaki}, Chapter I, page 6.

Now we can move on to approximation.

\begin{thm}\label{th:main}
Let $X$ and $Y$ be complex Banach spaces, $U \subset X$ open, $V \subset U$.  Suppose that the sequence $\{p_j\}$ is contained in a compact convex set $L \subset U$.  Let $W$ be the convex hull of $L \cup V$ and let $W'$ be the balanced convex hull of $L+V$. Suppose that for some $\rho > e$, $W+\rho W' \subset U$ and $f: U \rightarrow Y$ is holomorphic.  Then the Kergin series for $f$ converges to $f$ uniformly on compact subsets of V.
\end{thm}
\proof
First, we observe that if $T^k=(\mathbb{R}/\mathbb{Z})^k$
is the $k$-dimensional torus with Haar probability measure $dt$,
then
\begin{equation}\label{eq:homogeneity}
d^k f(a;v_0,\ldots, v_{k-1})=\int_{T^k} f(a+v_0 e^{2\pi i t_0}+ \ldots + v_{k-1} e^{2 \pi i t_{k-1}}) dt,
\end{equation}
provided each $v_j$ is small enough so that the right-hand side of \eqref{eq:homogeneity} is defined.
Let $s=(s_1,\ldots, s_k), a(s,x)=s_0 p_0+s_1 p_1+ \ldots + s_k  x, $ with $s_0$ again as in \eqref{eq:s_0}.  Plugging this into the error formula \eqref{eq:error} yields
\[
\left(\frac{k}{\rho}\right)^k \int_{S_k}  \int_{T^k} f \left(a(s,x)+ (x-p_0)\frac{\rho e^{2 \pi i t_0}}{k}+ \ldots +  (x-p_{k-1})\frac{\rho e^{2 \pi i t_{k-1}}}{k} \right)dtds.
\]
We have used homogeneity to factor out $\left(\frac{k}{\rho}\right)^k$.  Define
\[
c(s,t)=(x-p_0)\frac{\rho e^{2 \pi i t_0}}{k}+ \ldots + (x-p_{k-1})\frac{\rho e^{2 \pi i t_{k-1}}}{k}.
\]
Set $b(s,t)=a(s,x)+c(s,t)$.
If the interpolation points $p_j$ are in $L$ and $x \in V$, we can see that $a(s,x) \in W$ and $c(s,t) \in \rho W'$.  Thus, by hypothesis, we have that $b(s,t)$ takes values in $U$.

Now we restrict our attention to the case where $x$ is in a compact subset $K$ of $V$.
Applying Proposition \ref{pr:cpct}, plus the fact that the sum of two compact sets is compact, we see that the image of $b(s,t)$ is contained in a compact subset of $U$. Furthermore, this set is independent of $k$, and depends only on $K$ and $L$.
Let $M$ be the maximum of $f$ on this subset.  This leads to the inequality
\[
\left\| \left(\frac{k}{\rho}\right)^k \int_{S_k}  \int_{T^k} f(b(s,t))dtds \right\| \le
\left(\frac{k}{\rho}\right)^k M\vol (S_k).
\]
Recall that $\vol (S_k)=1/k!$, which, by Stirling's formula, is asymptotically equal to
\[
\left(\frac{e}{k}\right)^k \frac{1}{\sqrt{2\pi k}}.
\]
Substituting this all into the above estimate on the error of the $(k-1)^{st}$ degree Kergin polynomial yields
\[
\left(\frac{e}{\rho}\right)^k \frac{M}{\sqrt{2 \pi k}}.
\]
  Since $\rho>e,$ the error goes to zero. \proved

This allows $f$ to be unbounded if $V$ contains an open set, as
in the following example.
\begin{ex}
Let $B(r)=\{x \in X: \|x \|<r\}.$
Let $f:B(1) \rightarrow Y$ be holomorphic, and choose $r, r'>0$ such that
\begin{equation}\label{eq:rr'}
\frac{1-\max(r,r')}{r+r'}>e.
\end{equation}
Suppose that the closure of the set of all interpolation points $p_k$ is a compact subset of
$B(r').$
Then the corresponding Kergin series for $f$ converges to $f$ uniformly on compact subsets of $B(r).$
\end{ex}
\proof
Let $U=B(1),$ $V=B(r')$, and
\[
e < \rho < \frac{1-\max(r,r')}{r+r'}.
\]
Let $L$ be the convex hull of the closure of the interpolation points.  By Theorem \ref{th:main}, it suffices to check that $W+\rho W' \subset U$, where $W$ and $W'$ are as in the theorem.  Since $W$ is the convex hull of $L \cup V$, we have $W \subset B(\max(r,r'))$.  Furthermore, by the triangle inequality, we have that $W' \subset B(r+r')$.  Thus, if $x \in W+\rho W'$, then
\[
\|x\| < \max(r,r')+ \rho (r+r') \le \max(r,r') + \frac{1-\max(r,r')}{r+r'}(r+r')=1.
\]
In other words, $W+\rho W' \subset U$, as required.
\proved
Note that \eqref{eq:rr'} is satisfied if, for example,
\[
r=r'=\frac{1}{e+1}.
\]
Observe that here, $f$ may be unbounded even on a small ball, whereas in \cite{Filipsson}, $f$ must be bounded.  The price for such a convergence result is the stronger restriction on the interpolation points than the one found in \cite{Filipsson}.  Furthermore, Filipsson's convergence is uniform on balls, whereas Theorem \ref{th:main} only shows convergence on compact sets.  Now we give an example where the interpolation points are not in a compact set, and the Kergin series of an entire function diverges at the origin.
\begin{ex}
Let $f:l^1\rightarrow \mathbb{C}$ be defined by
\[
f(x)=\sum_{n=1}^\infty \left( n!\prod_{k=1}^n x_k \right).
\]
The function $f$ is entire, i.e., it is holomorphic on all of $l^1.$ Let $\{e_k\}$ be the standard basis for $l^1$ (0's everywhere except for a 1 in the $k^{th}$ position).  Using this basis as interpolation points yields a Kergin series that diverges at the origin.
\end{ex}
\proof
First, we show that $f$ is entire.
Let $x^0=(x_n^0) \in l^1.$
Choose $n_0$ so large that
\[
\sum_{n>n_0} |x_n^0| <\delta < e
\]
and let
\[
M>\prod_{n=1}^{n_0} |x_n^0|.
\]
We will show that the sum
\[
\sum_{n>n_0} n! \prod_{k=1}^n |x_k|
\]
converges uniformly near $x^0$, in fact on the set
\[
\left\{ x: \prod_{n=1}^{n_0}  |x_n| < M, \sum_{n>n_0} |x_n| < \delta \right\}.
\]
The modulus of the $(n_0+n)^{th}$ term in the sum defining $f$ is
\[
(n_0+n)!\prod_{k=1}^{n_0+n} |x_k| \le (n_0+n)! M \prod_{k=1}^n |x_{n_0+k}|
\]
\begin{equation}\label{eq:mean}
\le (n_0+n)! M \left(\frac{1}{n} \sum_{k=1}^n |x_{n_0+k}| \right)^n
\end{equation}
\begin{equation}\label{eq:Stirling}
\le (n_0+n)! M \left(\frac{\delta}{n}  \right)^n.
\end{equation}
In \eqref{eq:mean}, we used the fact that the geometric mean is less than the arithmetic mean.  Stirling's formula implies that \eqref{eq:Stirling} is asymptotically equal to
\[
\left( \frac{n_0+n}{e} \right)^{n_0+n} \sqrt{2 \pi (n_0+n)} M
\left( \frac{\delta}{n} \right)^n
\]
\[
=\left( \frac{n_0+n}{n} \right)^{n_0+n} \sqrt{2 \pi (n_0+n)} M
 \frac{\delta^n n^{n_0}}{e^{n+n_0}} =O \left(n^{n_0+1/2} \left(\frac{\delta}{e} \right)^n \right),
\]
because
\[
\left( \frac{n_0+n}{n} \right)^{n_0+n}= \left( 1+ \frac{n_0}{n} \right)^{n_0+n} \rightarrow e^{n_0}
\]
as $n \rightarrow \infty$.

Hence the series of $f$ converges uniformly near $x_0$ and so
f is entire.

Now we show that the Kergin series diverges at the origin.
Define $f_n=x_1 x_2 \cdots x_n$.  We check by induction that
\[
  d^k f_n (y;e_1, e_2, \ldots, e_k)=
    \left\{
      \begin{array}{ll}
        y_{k+1} \cdots y_n  & \text{ if } k<n, \\
        1                   & \text{ if } k=n, \\
        0                   & \text{ if } k>n.
      \end{array}
    \right.
\]

In particular, the left-hand side is non-negative when all $y_j \ge 0$.  Hence for such $y$
\[
(-1)^k d^k f(y;-e_1, -e_2, \ldots, -e_k) = \sum_{n=1}^\infty n! d^k f_n (y; e_1, \ldots , e_k) \ge  k!.
\]
Setting $y=s_0 e_1+ \ldots + s_k e_{k+1}$ and integrating over $S_k$ yields the absolute value of the $k+1^{st}$ term in the Kergin series at the origin, which must be at least $1$.  Summing, we have a divergent series.
\proved
Theorem \ref{th:main} requires that the interpolation sequence be contained in a compact subset of $U$.  Here, we set $U=l^1,$ $V=\{0\}.$  The only criterion not satisfied in the theorem is that $L$ must be compact.

\subsection*{Acknowledgments}
The author would like to thank L\'aszl\'o Lempert for bringing
Kergin interpolation to the author's attention, as well as for
his many useful suggestions.

\end{document}